\documentclass{amsart}
\usepackage{epsf}
\usepackage[all]{xy}
\usepackage{enumerate}
\usepackage{amssymb}
\usepackage{color}

\newcommand{\ie}{{\em i.e.\ }}

\newtheorem{theorem}{Theorem}
\newtheorem{lemma}{Lemma}
\newtheorem{proposition}{Proposition}
\newtheorem{proposition*}{Proposition}
\newtheorem{corollary}{Corollary}

\theoremstyle{definition}
\newtheorem{definition}{Definition}
\newtheorem{example}{Example}

\theoremstyle{remark}
\newtheorem*{remark}{Remark}

\numberwithin{equation}{section}

\newcommand{\internalcomment}[1]{}

\newcommand{\Z}{\mathbf{Z}}

\newcommand{\ko}{\: , \;}

\newcommand{\ol}{\overline}
\newcommand{\ul}{\underline}

\renewcommand{\tilde}[1]{\widetilde{#1}}

\newcommand{\ra}{\rightarrow}

\newcommand{\arr}[1]{\stackrel{#1}{\rightarrow}}

\newcommand{\opname}[1]{\operatorname{\mathsf{#1}}}

\renewcommand{\mod}{\opname{mod}\nolimits}
\newcommand{\Mod}{\opname{Mod}\nolimits}

\newcommand{\id}{\mathbf{1}}

\newcommand{\im}{\opname{im}\nolimits}
\renewcommand{\ker}{\opname{ker}\nolimits}

\newcommand{\op}[1]{\opname{#1}\nolimits}

\newcommand{\ca}{{\mathcal A}}
\newcommand{\cb}{{\mathcal B}}
\newcommand{\cc}{{\mathcal C}}
\newcommand{\cd}{{\mathcal D}}
\newcommand{\ce}{{\mathcal E}}

\newcommand{\ch}{{\mathcal H}}

\newcommand{\ct}{{\mathcal T}}
\newcommand{\cu}{{\mathcal U}}
\newcommand{\cv}{{\mathcal V}}

\newcommand{\cx}{{\mathcal X}}

\renewcommand{\phi}{\varphi}

\newcommand{\Hom}{\opname{Hom}}

\newcommand{\Ext}{\opname{Ext}}

\begin{document}
\title{Balance in stable categories}

\author{Pedro Nicol\'{a}s}

\address{Departamento de Matem\'{a}ticas Universidad de Murcia, Aptdo. 4021, 30100 Espinardo, Murcia, SPAIN}

\email{pedronz@um.es}

\author{Manuel Saor\'{\i}n}

\address{Departamento de Matem\'{a}ticas Universidad de Murcia, Aptdo. 4021, 30100 Espinardo, Murcia, SPAIN}

\email{msaorinc@um.es}
\thanks{The authors have been supported by the D.G.I. of the Spanish Ministry of Education
and the Fundaci\'{o}n `S\'{e}neca' of Murcia, with a part of FEDER
funds from the European Commission}

\begin{abstract}
We study when the stable category $\frac{\ca}{<\ct>}$ of an abelian category $\ca$ modulo a full additive subcategory $\ct$ is
balanced and, in case $\ct$ is functorially finite in $\ca$, we study a weak version of balance for $\frac{\ca}{<\ct>}$. Precise
necessary and sufficient conditions are given in case $\ct$ is either a Serre class or a class consisting of projective objects.
The results in this second case apply very neatly to (generalizations of) hereditary abelian categories.
\end{abstract}

\maketitle
\tableofcontents

\section{Introduction}\label{Introduction}
\bigskip

Given an additive category $\cc$ and a full additive subcategory $\ct$, the \emph{stable category} $\ul{\cc}:=\cc/{<\ct >}$ of $\cc$ with respect
to $\ct$ has the same objects as $\cc$, but $\ul{\cc}(X,Y)$ is the abelian group obtained from $\cc (X,Y)$ by making
zero those morphisms which factor through an object of $\ct$. These categories have played a fundamental role in Mathematics,
specially when $\cc =\ca$ is an abelian category and the subcategory $\ct$ has good homological properties. One  finds in
the literature early important applications of stable categories to Commutative Algebra and Algebraic Geometry \cite{AuslanderBridger1969} or Representation Theory of Artin algebras \cite{AuslanderReiten1974}, in this later case leading
to a milestone in the area, namely, the proof of the existence of almost split sequences in the category of finitely generated
modules over such algebras. In more recent times, a particular case has received much interest, namely,  when $\cc=\ca$ is
an abelian Frobenius category and $\ct$ is the subcategory of its projective (=injective) objects. That includes the case in which
$\ca =\cc\cb $ is the category of cochain complexes of an abelian category $\cb$, thus leading to the homotopy category $\ch\cb$
\cite[Chapter 1]{Happel1988}, \cite{Keller1996}, and the case in which $\ca =B\Mod$ is the category of modules over a block
of a  finite-group algebra. In this way, stable categories are a fundamental tool in modern Modular Representation Theory to tackle
long-standing open problems like the abelian defect group conjecture \cite{Rickard2001}, \cite{Rouquier2001}.

In full generality, stable categories are neither abelian nor
triangulated and, for a time,  there was a need to find a suitable
framework in which they could be inscribed. With sufficiently
general hypotheses, that framework is provided by one-sided
triangulated categories, in the terminology of
\cite{BeligianisMarmaridis1994}, or, equivalently, by suspended
categories, in the terminology of \cite{KellerVossieck1987}. As
shown by Beligianis and Marmaridis
\cite{BeligianisMarmaridis1994}, if $\ct$ is contravariantly
(resp. covariantly) finite in the abelian category $\ca$, then
$\ul{\ca}$ has a structure of left (resp. right) triangulated
category. In addition, if $\ct$ is functorially finite, \ie both
contravariantly and covariantly finite, then the corresponding
stable category satisfies all the axioms needed to be a
pretriangulated category. This latter class of categories includes
the abelian and triangulated ones and, also, the homotopy
categories of additive closed model categories in the sense of Quillen \cite[Chapter II]{BeligianisReiten2001}.

One of the common features of abelian and triangulated categories is that they are \emph{balanced}, \ie every morphism which is both a monomorphism and an epimorphism is necessarily an isomorphism. That property is no longer shared by pretriangulated categories, although monomorphisms and epimorphisms in them are easily detected using triangles. Indeed, the following lemma and its dual have an easy proof.

\begin{lemma} \label{mono in left-triangulated}
Let $\cx$ be a left triangulated category with loop functor $\Omega :\cx\ra\cx$. For a morphism $f:X\ra Y$ in $\cx$, the following assertions are
equivalent:
\begin{enumerate}[1)]
\item $f$ is a monomorphism.
\item There is a left triangle $\Omega Y\ra Z\arr{0}X\arr{f}Y$.
\end{enumerate}
\end{lemma}

In view of the above lemma, it makes sense to give the following:

\begin{definition}
Let $\cx$ be a left triangulated category with loop functor $\Omega :\cx\ra\cx$. A morphism $f:X\ra Y$ in $\cx$ will be said to be a \emph{strong monomorphism} if there is a left triangle $\Omega Y\ra 0\ra X\arr{f}Y$ in $\cx$. The dual notion on a right triangulated category is that of \emph{strong epimorphism}. Finally, if $\cx$ is a pretriangulated category we shall say that it is \emph{weakly balanced} in case that every morphism which is
both a strong mono and strong epi is an isomorphism.
\end{definition}

It is natural to expect that the class of balanced (resp. weakly
balanced) pretriangulated categories is an interesting one. The
moral goal of this paper is to show that (weak) balance is a very
restrictive condition on stable categories. We show that by
considering an arbitrary abelian category $\ca$ and making two
particular choices of the subcategory $\ct$, namely, when $\ct$ is
a Serre class or when $\ct$ consists of projective objects. The
organization of the paper goes as follows. In Section
\ref{Monomorphisms and strong monomorphisms in stable categories}
we characterize monomorphism and strong monomorphisms in
$\ul{\ca}$, for any full additive subcategory $\ct$ closed under
direct summands. In Section \ref{When T is a Serre class} we show
that if $\ct$ is a covariantly finite Serre class in $\ca$, then
$\ul{\ca}$ is (weakly) balanced if and only if $\ct$ is a direct
summand  of $\ca$ as an additive category (Corollary \ref{balanced
in Serre case}). In the later sections of the paper $\ct$ is
supposed to consist of projective objects, in which case the
balance of $\ul{\ca}$ is characterized in Section \ref{Balance
when T consists of projective objects} (Theorem
\ref{characterization balanced}) and the weak balance in Section
\ref{Weak balance when T consists of projective objects} (Theorem
\ref{characterization weakly-balanced}). A byproduct is that if
$\ct$ consists of projective-injective objects, then $\ul{\ca}$ is
always balanced. When $\ca=\ch$ has the property that subobjects
of projective objects are projective (\emph{e.g.}\ if $\ch$ is a
hereditary abelian category), these results give that if $\ct$ is
covariantly finite then $\ul{\ch}$ is (weakly) balanced if and
only if $^\perp\ct$ is closed under taking subobjects (Proposition
\ref{hereditary weakly balanced}). The final section of the paper
gives some examples in categories of modules showing, in
particular, that stable balanced categories need not be either
abelian or triangulated and, also, that weakly balanced stable
categories need not be balanced.

All throughout the paper $\ca$ will be an abelian category,  $\ct$ will be  a full additive subcategory closed under direct summands, ${<\ct>}$ will be the ideal of $\ca$ formed by the morphisms that factors through an object of $\ct$ and $\ul{\ca}:=\ca/{<\ct>}$ the associated quotient category.  The
hypothesis of $\ct$ being closed under direct summands is not strictly necessary for the ideals ${<\ct>}$ and ${<\op{add}(\ct )>}$ coincide, where $\op{add}(?)$ denotes the closure under finite direct sums and direct summands.  But it simplifies the statements and we will assume it in the sequel.

The reader is referred to \cite{RadaSaorinDelValle2000} for the concepts of (pre)envelopes, (pre)covers and contravariantly (resp. covariantly) finite subcategories that we shall frequently use in the paper.  By \cite[Theorem 3.1]{BeligianisMarmaridis1994}, when $\ct$ is contravariantly (resp. covariantly) finite in $\ca$,  the stable category $\ul{\ca}:=\ca/{<\ct>}$ has a canonical structure of left (resp. right) triangulated category, where the loop (resp.
suspension) functor $\Omega :\ul{\ca}\ra\ul{\ca}$ (resp. $\Sigma :\ul{\ca}\ra\ul{\ca}$) is defined on objects as the kernel (resp. cokernel) of a $\ct$-precover (resp. $\ct$-preenvelope). When $\ct$ is functorially finite in $\ca$, the category $\ul{\ca}$ becomes pretriangulated. The reader is referred to \cite{KellerVossieck1987}, \cite{BeligianisMarmaridis1994} and \cite{BeligianisReiten2001} for the definitions and terminology concerning one-sided
triangulated (or suspended) categories and pretriangulated categories.

The terminology on general categories used in this paper is standard and can be found in common textbooks on the subject, like \cite{Pareigis1973}. For the sake of intuition and simplification of proofs, we shall allow ourselves some abuse of terminology concerning abelian categories and shall use expressions of module theory, like ``canonical inclusion $\im(f)\hookrightarrow Y$'' for a morphism $f:X\ra Y$, meaning the (unique up to isomoprhism) monomorphism appearing in the epi-mono factorization of $f$. We hope that our choice will not stop the reader from understanding statements and proofs. For the final section, the terminology concerning rings can be found in \cite{AndersonFuller1992} and \cite{Rotman1979} while that concerning finite dimensional algebras can be found in \cite{Ringel1984} and \cite{AuslanderReitenSmalo1995}.
\bigskip

\textbf{Acknowledgements:} The authors thank Kent Fuller for telling them about the main result of reference 8 and
to Bernhard Keller for calling their attention on the result of Verdier mentioned in the
final remark of Section \ref{Weak balance when T consists of projective objects}.
\bigskip

\section{Monomorphisms and strong monomorphisms in stable categories}\label{Monomorphisms and strong monomorphisms in stable categories}
\bigskip

In this section we undertake the identification of monomorphisms
and strong monomorphisms in $\ul{\ca}$. All throughout the paper,
the image of a morphism $f$ under the canonical projection
$\ca\longrightarrow\ul{\ca}=\ca/<\ct >$ will be denoted by
$\ol{f}$.

\begin{proposition}\label{monos}
For a morphism $f:X\ra Y$ in $\ca$ the following assertions are equivalent:
\begin{enumerate}[1)]
\item $\ol{f}$ is a monomorphism in $\ul{\ca}$.
\item For every $p\in\ca(T,Y)$ with $T\in\ct$, the parallel to $p$ in the pullback of $f$ and $p$ factors through an object of $\ct$.
\end{enumerate}

When, in addition, $\ct$ is contravariantly finite in $\ca$ the following assertions are also equivalent:

\begin{enumerate}[3)]
\item For some $\ct$-precover $p_{Y}:T_{Y}\ra Y$ the parallel to $p_{Y}$ in the pullback of $f$ and $p_{Y}$ factors through an object of $\ct$.
\end{enumerate}
\begin{enumerate}[4)]
\item In the left triangulated category $\ul{\ca}$ there is a left triangle of the form
\[\Omega Y\ra Z\arr{\ol{0}} X\arr{\ol{f}} Y
\]
\end{enumerate}
\end{proposition}
\begin{proof}
$1)\Rightarrow 2)$ From the cartesian square
\[\xymatrix{Z\ar[r]^{g}\ar[d]^{h} & X\ar[d]^{f} \\
T\ar[r]_{p} & Y
}
\]
we have $\ol{f}\ol{g}=\ol{0}$ and so $\ol{g}=\ol{0}$, \ie $g$ factors through an object of $\ct$.

$2)\Rightarrow 1)$ Assume $\ol{f}\ol{u}=\ol{0}$ for some morphism $u:U\ra X$. Then we have a commutative square
\[\xymatrix{U\ar[r]^{u}\ar[d]_{v} & X\ar[d]^{f} \\
T\ar[r]_{p} & Y
}
\]
By doing the pullback of $f$ and $p$ we get the following commutative diagram
\[\xymatrix{U\ar@{.>}[dr]_{w}\ar@/^/[drr]^{u}\ar@/_/[ddr]_{v} & & \\
& Z\ar[r]^{g}\ar[d]_{h} & X\ar[d]^{f} \\
& T\ar[r]_{p} & Y
}
\]
Since, by hypothesis, $\ol{g}=\ol{0}$, and we have $\ol{u}=\ol{g}\ol{w}$, we conclude $\ol{u}=\ol{0}$.

$1)\Leftrightarrow 4)$ (see Lemma \ref{mono in left-triangulated}) and $2)\Rightarrow 3)$ are clear. As for $3)\Rightarrow 4)$, use the canonical construction of left triangles (cf. \cite{BeligianisMarmaridis1994} or \cite[Chapter II.1]{BeligianisReiten2001}).
\end{proof}

The following result identifies strong monomorphisms in stable categories.

\begin{proposition}\label{left zero}
Assume that $\ct$ is contravariantly finite in $\ca$. The following assertions are equivalent for a morphism $f\in\ca(X,Y)$:
\begin{enumerate}[1)]
\item $\ol{f}$ is a strong monomorphism in the left triangulated category $\ul{\ca}$.
\item For every (resp. some) $\ct$-precover $p_{Y}:T_{Y}\ra Y$ the parallel to $p_{Y}$ in the pullback of $f$ and $p_{Y}$ is a $\ct$-precover of $X$.
\item For every (resp. some) $\ct$-precover $p_{Y}:T_{Y}\ra Y$, in the cartesian square
\[\xymatrix{T\ar[r]^g\ar[d]_{h} & X\ar[d]^{f} \\
T_{Y}\ar[r]_{p_{Y}} & Y
}
\]
we have $T\in\ct$.
\end{enumerate}
\end{proposition}
\begin{proof}
$2)\Rightarrow 3)$ is clear and $3)\Leftrightarrow 1)$ follows from the construction of left triangles in $\ul{\ca}$ \cite[Chapter II.1]{BeligianisReiten2001}.

$3)\Rightarrow 2)$ Consider the cartesian square of $3)$ and a morphism $g':T'\ra X$ with $T'\in\ct$. Since $p_{Y}$ is a $\ct$-precover,
there exists a morphism $h':T'\ra T_{Y}$ such that $fg'=p_{Y}h'$. Now, by the universal property of the cartesian square, there exists
a morphism $t:T'\ra T$ such that $gt=g'$.
\end{proof}

\begin{remark}\label{observaciones}
\begin{enumerate}[1)]
\item Applying the duality principle, we obtain corresponding results, dual to  the above,  when considering the right triangulated structure in $\ul{\ca}$ defined by a covariantly finite additive full subcategory $\ct$ of $\ca$ closed under direct summands. We leave the statements for the reader.
\item Under the hypotheses of Proposition ~\ref{left zero}, if $f:X\ra Y$ is a morphism in $\ca$ such that $\ol{f}$ is a strong monomorphism  in
$\ul{\ca}$, then $\ol{f}=\ol{g}$, for some morphism $g:X'\ra Y$ such that $\ker(g)\in\ct$. Indeed, take $g=\scriptsize{\left[\begin{array}{cc}f&
p_Y\end{array}\right]}:X\oplus T_Y\ra Y$, where $p_Y:T_Y\ra Y$ is a $\ct$-precover.
\item If, in the situation of $(2)$, the subcategory $\ct$ consists of injective objects of $\ca$, then $g$ can be chosen to be a monomorphism in $\ca$. Indeed, according to $(2)$, we can assume from the beginning that $\ker(f)\in\ct$. Then the monomorphism $\ker(f)\ra X$ is a section and we can write
$f=\scriptsize{\left[\begin{array}{cc}0&f_{|X'}\end{array}\right]}:\ker(f)\oplus X'\ra Y$, for some complement $X'$ of $\ker(f)$ in $X$. Now choose
$g=f_{|X'}$
\end{enumerate}
\end{remark}
\bigskip

\section{When $\ct$ is a Serre class}\label{When T is a Serre class}
\bigskip

 Recall that, given an abelian category $\ca$, a full subcategory $\ct\subseteq\ca$ is called a \emph{Serre class} if, given an exact sequence $0\ra X\ra Y\ra
Z\ra 0$ in $\ca$, the object $Y$ belongs to $\ct$ if and only if $X$ and $Z$ belong to $\ct$. All throughout this section, we assume that $\ct$ is a Serre class in $\ca$. The selfdual condition of Serre classes allows to give a dual result of any one that we shall give here.

Recall also that a pair $(\cu ,\cv )$ of full subcategories is a \emph{torsion theory} in $\ca$ in case $\cu^\perp:=\{X\in\ca\mid \ca (U,X)=0\text{ for all }U\in\cu\} =\cv\ko \cu =^\perp\cv:=\{X\in\ca\mid\ca (X,V)=0\text{ for all }V\in\cv\}$ and the inclusion functor $i:\cu\hookrightarrow\ca$ has a right
adjoint. The concept is selfdual since in that case the inclusion functor $j:\cv\hookrightarrow\ca$ has a left adjoint.  The following is an auxiliary lemma.

\begin{lemma} \label{(co- contra-)variantly finite Serre classes}
Let $\ca$ be an abelian category and $\cu$ be a full subcategory closed under taking quotients and extensions. The following assertions are equivalent:
\begin{enumerate}[1)]
\item $\cu$ is contravariantly finite in $\ca$.
\item The inclusion functor $i:\cu\hookrightarrow\ca$ has a right adjoint.
\item $(\cu ,\cu^\perp )$ is a  torsion theory in $\ca$.
\end{enumerate}
\end{lemma}
\begin{proof}
$3)\Rightarrow 2)$ and $2)\Rightarrow 1)$ are clear. We just need to prove $1)\Rightarrow 3)$. If $p_X:U_X\ra X$ is a $\cu$-precover of $X$, then $\im(p_X)=:u(X)$ belongs to $\cu$ and one easily shows that the inclusion morphism $j_X:u(X)\ra X$ is a $\cu$-cover, whence uniquely determined up to isomorphism. We leave to the reader the easy verification that the assignment $X\mapsto u(X)$ extends to a functor $u:\ca\ra\cu$ which is right adjoint to
the inclusion.

Now $\cu^\perp =\{F\in\ca\mid u(F)=0\}$ and we just need to prove that $\cu =^\perp (\cu^\perp )$ or, what is enough, that $u(X/u(X))=0$ for all $X\in\ca$. But this is clear for $u(X/u(X))={U}/{u(X)}$ for some subobject $U$ of $X$ containing $u(X)$ and fitting in a short exact sequence $0\rightarrow u(X)\ra U\ra
u(X/u(X))\rightarrow 0$. Since the two outer terms belong to $\cu$ we also have $U\in\cu$. But then the canonical inclusion $U\hookrightarrow X$ factors through $j_X:u(X)\ra X$, which implies that $U\subseteq u(X)$ and, hence, we get $u(X/u(X))=0$ as desired.
\end{proof}

Our next result identifies the monomorphisms and epimorphisms in $\ul{\ca}$.

\begin{proposition} \label{triangles for Serre classes}
The following statements hold for a morphism $f:X\ra Y$ in $\ca$:
\begin{enumerate}[1)]
\item $\ol{f}$ is a monomorphism in $\ul{\ca}$ if and only if $\ker(f)$ belongs to $\ct$. In case $\ct$ is contravariantly
finite in $\ca$, that happens if and only if $\ol{f}$ is a strong monomorphism in $\ul{\ca}$.
\item[1')] $\ol{f}$ is an epimorphism in $\ul{\ca}$ if and only if $\op{cok}(f)$ belongs to $\ct$. In case $\ct$ is covariantly finite in $\ca$, that happens if and
only if $\ol{f}$ is a strong epimorphism in $\ul{\ca}$.
\item If $X\in\ct^\perp$, then $\ol{f}$ is an isomorphism in $\ul{\ca}$ if and only if $f$ is a section in $\ca$ with $\op{cok}(f)\in\ct$.
\end{enumerate}
\end{proposition}
\begin{proof}
$1')$ follows from $1)$ by duality. So we prove assertion $1)$. If $\ol{f}$ is a monomorphism in $\ul{\ca}$, then we know that we
have a factorization as follows:
\[f^k: \ker(f   )\arr{j}T\arr{v} X
\]
with $T\in\ct$ (notice that $j$ is a monomorphism). Since $\ct$ is closed for subobjects, we conclude that $\ker(f)\in\ct$.

Suppose now that $\ker(f)\in\ct$. According to Proposition ~\ref{monos}, we need to prove that, for every morphism $h:T\ra Y$
with $T\in\ct$, the parallel to $h$ in the pullback of $f$ and $h$ factors through an object of $\ct$. But here we have even more,
for if $Z$ is the upper left corner of that pullback, we have an exact sequence $0\ra \ker(f)\ra Z\arr{h}T$. The fact that $\ct$ is
a Serre class implies that $Z\in\ct$. This same argument also proves, with the help of Proposition ~\ref{left zero}, that if
$\ct$ is covariantly finite, then $\ol{f}$ is a strong monomorphism in $\ul{\ca}$. That ends the proof of assertion $1)$.

We next prove (the 'only if' part of) assertion $2)$. If $X\in\ct^\perp$ and $\ol{f}$ is an isomorphism in $\ul{\ca}$, then
assertions $1)$ and $1')$ imply that $\ker(f), \op{cok}(f)\in\ct$. But then $\ker(f)\in\ct\cap\ct^\perp$, which implies $\ker(f)=0$. On the other hand, since $\ol{f}$ is an isomorphism there exists a morphism $g:Y\ra X$ such that $\id_X-gf$ factors through an object of $\ct$. That is, we
have a factorization $\id_X-gf:X\ra T\ra X$. But, since $X\in\ct^\perp$,  the second morphism of that factorization is
zero, so that $\id_X-gf=0$  and, hence,  $f$ is a section in $\ca$ with $\op{cok}(f)\in\ct$.
\end{proof}

\begin{definition}
Let $\cc$ be any additive category and $\cd$, $\ce$ two full subcategories. We shall say that $\cc$ is the \emph{direct sum} of $\cd$ and $\ce$, written $\cc =\cd\oplus\ce$, when the following two conditions hold:
\begin{enumerate}[i)]
\item Every object of $\cc$ decomposes as a direct sum of an object of $\cd$ and an object of $\ce$.
\item $\cc(D,E)=0=\cc(E,D)$, for all objects $D\in\cd$ and $E\in\ce$
\end{enumerate}
\end{definition}

The main result of this section follows now in an easy way:

\begin{corollary} \label{balanced in Serre case}
Let $\ca$ be an abelian category and $\ct$ be a Serre class in $\ca$. The following assertions are equivalent:
\begin{enumerate}[1)]
\item $\ct$ is contravariantly (resp. covariantly) finite in $\ca$ and $\ul{\ca}$ is balanced.
\item $\ct$ is functorially finite in $\ca$ and $\ul{\ca}$ is weakly balanced.
\item $\ca=\ct\oplus\ct^\perp$.
\end{enumerate}
\end{corollary}
\begin{proof}
Since assertion $2)$ is selfdual we just consider the contravariantly finite version of assertion $1)$.

$2)\Rightarrow 1)$ follows from Proposition ~\ref{triangles for Serre classes}.

$3)\Rightarrow 2)$ It is easy to see that $\ct$ is functorially finite in $\ca$, that we have an equivalence $\ca/{<\ct>}\simeq
\ct^{\perp}$ of additive categories and that $\ct^{\perp}$ is abelian (and so $\ca/{<\ct>}$ will be balanced).

$1)\Rightarrow 3)$  We first prove that $\ca (F,T)=0$, for all $F\in\ct^\perp$ and $T\in\ct$ or, equivalently, that
$\ct^\perp\subseteq ^\perp\ct$. Indeed, let $f:F\ra T$ be a morphism. Replacing $T$ by $\im(f)$ if necessary, we can
assume that $f$ is an epimorphism in $\ca$ and we need to prove that $T=0$. For that we consider the monomorphism $f^k:\ker(f)\ra F$. By Proposition ~\ref{triangles for Serre classes}, we know that $\ol{f^k}$ is both a monomorphism and an epimorphism in $\ul{\ca}$. Since this is a balanced category, we conclude that $\ol{f^k}$ is an isomorphism. But then Proposition ~\ref{triangles for Serre classes} gives that $f^k$ is a section (with cokernel in
$\ct$). Then $f$ is a retraction, but any section for it must be zero because $T\in\ct$ and $F\in\ct^\perp$. Therefore $T=0$ as desired.

On the other hand, since $(\ct ,\ct^\perp )$ is a torsion theory (cf. Lemma \ref{(co- contra-)variantly finite Serre classes}) we can consider,  for every object $X\in\ca$, the canonical exact sequence $0\ra t(X)\ra X\arr{p_X}X/t(X)\ra 0$. According to Proposition \ref{triangles for Serre classes},  we know that
$\ol{p}_X$ is both a monomorphism and an epimorphism, whence an isomorphism,  in $\ul{\ca}$. Since $X/t(X)\in\ct^\perp\subseteq
^\perp\ct$ the dual of Proposition \ref{triangles for Serre classes}(2) says that $p_X$ is a retraction in $\ca$. Then we have
an isomorphism $X\cong t(X)\oplus\frac{X}{t(X)}$ and the proof is finished.
\end{proof}
\bigskip

\section{Balance when $\ct$ consists of projective objects}\label{Balance when T consists of projective objects}
\bigskip

Throughout the rest of the paper, unless explicitly said
otherwise,   the full subcategory $\ct$ of $\ca$ consists of
projective objects. Notice that we do not assume that $\ca$ has
enough projectives. Dual results of ours can be obtained by
assuming instead that $\ct$ consists of injective objects. We
leave their statement for the reader.

Our next observation is that every epimorphism in $\ul{\ca}$ can be represented by an epimorphism of $\ca$, something that was already noticed  by Auslander and Bridger in the (classical) projectively stable category of a module category (cf. \cite{AuslanderBridger1969} and, for related
questions, see also \cite{Kato2005}).

\begin{lemma}\label{reduction to epis}
Let $f\in\ca(X,Y)$ be such that $\ol{f}$ is an epimorphism in $\ul{\ca}$. Then the canonical monomorphism $\im(f)\ra Y$ yields
a retraction in $\ul{\ca}$ and there exists an epimorphism $f'\in\ca(X',Y)$ such that $\ol{f'}\cong\ol{f}$ in $\ul{\ca}$.
\end{lemma}
\begin{proof}
Let $f^c$ be the cokernel of $f$ in $\ca$. From $f^cf=0$ and the fact that $\ol{f}$ is an epimorphism we deduce that $\ol{f^c}=\ol{0}$. Therefore,
there exists a factorization as follows:
\[f^c:Y\arr{h}T\arr{p}\op{cok}(f)
\]
with $T\in\ct$.
%Observe that since $f^c$ is an epimorphism, then $p$ is an epimorphism and thus $\op{cok}(f)$ is always an epimorphic image of an object of $\ct$.

Now, since $T$ is projective there exists a morphism $g:T\ra Y$
such that $f^cg=p$, and so $f^c(gh-\id_{Y})=0$. Hence $gh-\id_{Y}$
factors through the canonical morphism $\im(f)\ra Y$, which is
the kernel of $f^c$. Then $\im(f)\ra Y$ yieds a retraction in
$\ul{\ca}$. This proves the first assertion.

Of course the morphism
$f':=\scriptsize{\left[\begin{array}{cc}f&g\end{array}\right]}:
X\oplus T\ra Y$ is a morphism such that  $\ol{f}'=\ol{f}$ in
$\ul{\ca}$. Let us see that $f'$ is an epimorphism in $\ca$. If we
have $\varphi f'=0$, then $\varphi f=0$ and $\varphi g=0$. Hence,
there exists a morphism $u$ such that $uf^c=\varphi$, which
implies $uf^cg=0$, \ie $up=0$. Since $p$ is an epimorphism we get
$u=0$ and so $\varphi=0$.
\end{proof}

Since our first goal is to characterize when $\ul{\ca}$ is
balanced, the above observation suggests to characterize those
epimorphisms in $\ca$ which become epimorphisms, monomorphisms or
isomorphisms when passing to $\ul{\ca}$. Our next results go in
that direction.

\begin{proposition}\label{epimorphisms which become stable-epis}
Let $f\in\ca(X,Y)$ be an epimorphism in $\ca$. The following assertions are equivalent:
\begin{enumerate}[1)]
\item $\ol{f}$ is an epimorphism in $\ul{\ca}$.
\item For every morphism $h:X\ra T$ with $T\in\ct$, there is a morphism $\tilde{h}:X\ra h(\ker(f))$ such that $h$ and $\tilde{h}$ coincide
on $\ker(f)$.

When $\ct$ is also covariantly finite, the above conditions are also equivalent to:

\item For every (resp. some) $\ct$-preenvelope $\mu^X:X\ra T^X$, there is a morphism $\tilde{\mu}:X\ra \mu^X(\ker(f))$ such that $\mu^X$ and
$\tilde{\mu}$ coincide on $\ker(f)$.
\end{enumerate}
\end{proposition}
\begin{proof}
$1)\Rightarrow 2)$ According to the dual of Proposition ~\ref{monos}, if $\ol{f}$ is an epimorphism then, for every
morphism $h:X\ra T$ to an object  $T\in \ct$, the parallel to $h$ in the pushout of $f$ and $h$ factors through an
object of $\ct$. We consider that pushout:
\[\xymatrix{X\ar[r]^{f}\ar[d]_{h} & Y\ar[d]^{v} \\
 T\ar[r]_{u} & Z
}\] so that we have a factorization $v:Y\arr{v_1}T'\arr{v_2}Z$,
with $T'\in\ct$. Since $u$ is an epimorphism in $\ca$ and $T'$ is
projective, we get a morphism $\varphi :T'\ra T$ such that
$u\varphi =v_2$. Then one gets $u\varphi v_1=v$. Hence
$u(h-\varphi v_1 f)=uh-vf=0$. From that we get a morphism
$\tilde{h} :X\ra \ker(u)$ such that $h-\varphi v_1
f=u^{k}\tilde{h}$, where $u^k$ is the kernel of $u$, and thus it
follows easily that $h$ and $\tilde{h}$ coincide on $\ker(f)$. We
only have to notice that, since $f$ is an epimorphism in $\ca$, a
classical construction of pushouts in abelian categories gives
that $u$ is the cokernel map of the composition $\ker(f)\arr{f^k}
X\arr{h} T$. Then $u^k:\ker(u)\ra T$ gets identified with the
canonical inclusion $h(\ker(f))\hookrightarrow T$.

$2)\Rightarrow 1)$ Let $g:Y\ra Z$ be a morphism in $\ca$ such that $\ol{g}\ol{f}=\ol{0}$ in $\ul{\ca}$. Then we
have a factorization $gf:X\arr{h}T\arr{p}Z$, where $T\in\ct$. The hypothesis says that we have a morphism
$\tilde{h}:X\ra h(\ker(f))$ such that $j\tilde{h}f^k=hf^k$, where $j:h(\ker(f))\ra T$ is the monomorphism of the canonical epi-mono factorization:
\[\xymatrix{\ker(f)\ar[rr]^{hf^k}\ar[dr]_{q} && T \\
& h(\ker(f))\ar[ur]_{j} &
}
\]
Then $(h-j\tilde{h})f^k=0$, which implies that $h-j\tilde{h}$
factors through the cokernel of $f^k$, \ie through $f$. Then we
have a morphism $\xi :Y\ra T$ such that $h-j\tilde{h}=\xi f$ and
we get $gf=ph=p(j\tilde{h}+\xi f)$. Since $pjq=phf^k=gff^k=0$ we
get that $pj=0$ and so $gf=p\xi f$. Being $f$ an epimorphism, we
conclude that $g=p\xi$,  and hence  $\ol{g}=\ol{0}$.

Suppose now that $\ct$ is covariantly finite.  We only have to prove (the weak version of) $3)\Rightarrow 2)$. Let us
fix a $\ct$-preenvelope $\mu^X:X\ra T^X$ for which condition $3)$ holds. If now $h:X\ra T$ is a morphism in
$\ca$, with $T\in\ct$, then there is a factorization $h:X\arr{\mu^X}T^X\arr{g}T$.
Then $g$ induces a morphism $\tilde{g}:\mu^X(\ker(f))\ra h(\ker(f))$. The hypothesis
says that we have a morphism $\tilde{\mu}:X\ra\mu^X(\ker(f))$ such that $\mu^X$ and
$\tilde{\mu}$ coincide on $\ker(f)$. Put now $\tilde{h}=:\tilde{g}\tilde{\mu}:X\ra h(\ker(f))$
and one easily checks that $h$ and $\tilde{h}$ coincide on $\ker(f)$.
\end{proof}

\begin{proposition} \label{epimorphisms which become stable-monos}
Let $f:X\ra Y$ be an epimorphism in $\ca$. Then $\ol{f}$ is a
monomorphism in $\ul{\ca}$ if and only if the canonical
monomorphism $f^k: \ker(f)\ra X$ factors through an object of
$\ct$.
\end{proposition}
\begin{proof}
If $\ol{f}$ is a monomorphism then, since
$\ol{f}\ol{f^k}=\ol{0}$, we have $\ol{f^k}=\ol{0}$, \ie $f^k$
factors through $\ct$. On the other hand, suppose that the
canonical monomorphism $f^k:\ker(f)\ra X$ factors in the form
$f^k:\ker(f)\arr{u}T\arr{q}X$, where $T\in\ct$. According to
Proposition \ref{monos}, it is enough to prove that, for every
morphism $p:T'\ra Y$  from an object $T'\in\ct$,  the parallel to
$p$ in the pullback of $f$ and $p$ is a morphism which factors
through an object of $\ct$. We then consider that pullback:
\[\xymatrix{ Z\ar[r]^{r}\ar[d]_{g} & X\ar[d]^{f} \\
T'\ar[r]_{p} & Y }
\]
Since $f$ is epi it follows that $g$ is also epi, and we have an
exact sequence
\[0\ra \ker(f)\ra Z\arr{g}T'\ra 0
\]
in $\ca$. Since $T'$ is projective this sequence splits and we can
identify $Z=\ker(f)\oplus T'$,
$g=\scriptsize{\left[\begin{array}{cc}0&\id\end{array}\right]}:\ker(f)\oplus
T'\ra T'$ and
$r=\scriptsize{\left[\begin{array}{cc}f^k&\xi\end{array}\right]}:\ker(f)\oplus
T'\ra X$, where $\xi:T'\ra X$ is a morphism such that $f\xi=p$. But then we have the following factorization of $r$:
\[\xymatrix{ \ker(f)\oplus T'\ar[rr]^{\scriptsize{\left[\begin{array}{cc}u&0\\0&\id\end{array}\right]}} &&
T\oplus
T'\ar[rr]^{\scriptsize{\left[\begin{array}{cc}q&\xi\end{array}\right]}}
&& X }
\]
Therefore $r$ factors through the object $T\oplus T'\in\ct$.
\end{proof}

\begin{corollary}\label{epimorphisms which become stable-monepis}
Let $f:X\ra Y$ be an epimorphism in $\ca$. The following
assertions are equivalent:
\begin{enumerate}[1)]
 \item $\ol{f}$ is both a monomorphism and an epimorphism in
$\ul{\ca}$ \item The canonical monomorphism $f^k:\ker(f)\ra X$
factors through an object of $\ct$ and, for every morphism $h:X\ra
T$ to an object of $\ct$, there is a morphism $\tilde{h}:X\ra
h(\ker(f))$ such that $h$ and $\tilde{h}$ coincide on $\ker(f)$.
\end{enumerate}
\end{corollary}

The following result finishes this series of preliminaries leading
to the characterization of when $\ul{\ca}$ is balanced.

\begin{proposition} \label{epimorphisms which become stable-isos}
 Let $f:X\ra Y$ be an
epimorphism in $\ca$. The following assertions are equivalent:
\begin{enumerate}[1)]
\item $\ol{f}$ is an isomorphism in $\ul{\ca}$. \item $f$ is a
retraction in $\ca$ with   $\ker(f)\in\ct$ \item $f$ is a
retraction in $\ca$ whose kernel map $f^k:\ker(f)\ra X$
factors through an object of $\ct$.
\end{enumerate}
\end{proposition}
\begin{proof}
$2)\Rightarrow 1)$ and $2)\Rightarrow 3)$ are clear.

$3)\Rightarrow 2)$ If 3) holds then $\ker(f)$ is a direct summand of an object of $\ct$, and our assumption on $\ct$ implies
that $\ker(f)\in\ct$.

$1)\Rightarrow 3)$ Since $\ol{f}$ is an isomorphism it is also a monomorphism and an epimorphism in $\ul{\ca}$, so that Corollary \ref{epimorphisms which become stable-monepis} applies. Let us choose $g:Y\ra X$ such that $\ol{g}=\ol{f}^{-1}$. Then $\id_X-gf$ factors in the form
$X\arr{h}T'\arr{v}X$, with $T'\in\ct$. Applying Corollary \ref{epimorphisms which become stable-monepis}, we get a morphism $\tilde{h}:X\ra h(\ker(f))$ such that $j\tilde{h} f^k=h f^k$, where $j:h(\ker(f))\hookrightarrow T'$ is the canonical inclusion (\ie the monomorphism in the epi-mono factorization of $h f^k$).

Notice now that the restriction of $v h=\id_X-g f$ to $\ker(f)$ is the identity map, in particular $(vh)(\ker(f))=\ker(f)$. If we denote by $\hat{v}$ the  morphism
$h(\ker(f))\ra (vh)(\ker(f))=\ker(f)$ induced by $v$ then, by definition, we have $f^k\hat{v}=v j$. But then
$f^k\hat{v}\tilde{h} f^k= vj\tilde{h} f^k=v h f^k=(\id_X-g f)f^k=f^k$.  Since $f^k$ is a monomorphism we then get that $\hat{v}\tilde{h} f^k$ is the identity on
$\ker(f)$. Therefore $f^k$ is a section which, according to Corollary \ref{epimorphisms which become stable-monepis}, factors through an
object of $\ct$.
\end{proof}

We come now to the main result of this section.

\begin{theorem}\label{characterization balanced}
Let $\ca$ be an abelian category and  $\ct$ be a full subcategory consisting of projective objects which is closed under taking
finite direct sums and direct summands. The following assertions are equivalent:
\begin{enumerate}[1)]
\item $\ul{\ca}$ is balanced.
\item If $\mu:K\ra X$ is a non-split monomorphism which  factors through an object of $\ct$,  then there exists a morphism $h:X\ra T'$, with $T'\in\ct$,
such that no morphism $\tilde{h}:X\ra (h\mu)(K)$ coincides with $h$ on $K$.
\item If $\mu:T\ra X$ is a non-split monomorphism with $T\in\ct$, then there exists a morphism $h:X\ra T'$, with $T'\in\ct$, such that no morphism $\tilde{h}:X\ra (h\mu)(T)$ coincides with $h$ on $T$.
\item If $f:X\ra Y$ is an epimorphism satisfying conditions $i)$ and $ii)$ below, then it is a retraction:
\begin{enumerate}
\item[i)] $f^k:\ker(f)\ra X$ factors through an object of $\ct$.
\item[ii)] For every $h:X\ra T\in\ct$ the canonical epimorphism $\ker(f)\twoheadrightarrow h(\ker(f))$ factors through $f^k:\ker(f)\ra X$.
\end{enumerate}
\item For every non-split epimorphism $f:X\ra Y$ with
$\ker(f)\in\ct$, there exists a morphism $h:X\ra T'$,  with
$T'\in\ct$,  such that $\ker(f)+\ker(h-gf)$ is strictly contained
in $X$ for all morphisms $g:Y\ra T'$.
\end{enumerate}
\end{theorem}
\begin{proof}
$4)\Leftrightarrow 2)$ and $2)\Rightarrow 3)$ are clear.

$1)\Leftrightarrow 4)$ is a consequence of Corollary \ref{epimorphisms which become stable-monepis} and Proposition
\ref{epimorphisms which become stable-isos}, bearing in mind that
every epimorphism of $\ul{\ca}$ can be represented by an
epimorphism of $\ca$.

$3)\Rightarrow 2)$ Let $\mu:K\ra X$ be a non-split monomorphism of $\ca$ which admits  a factorization
\[\mu:K\arr{u}T\arr{v}X
\]
with $T\in\ct$. Consider the following pushout diagram $(*)$ with exact rows
\[\xymatrix{0\ar[r] & K\ar[r]^{\mu}\ar[d]^{u} & X\ar[r]\ar[d]^{u'} & Y\ar[r]\ar@{=}[d] & 0 \\
0\ar[r] & T\ar[r]^{\mu'} & X'\ar[r] & Y\ar[r] & 0 }
\]
We work separately the case in which $\mu '$ does not split and
the case in which it  splits. In the first case, assertion $3)$ gives a morphism $h':X'\ra T'$, with $T'\in\ct$,  such that no
morphism $\tilde{h}':X'\ra h'\mu'(T)$ coincides with
$h'$ on $T$. We claim that $h=:h'u'$ satisfies the property needed
in assertion $2)$. Suppose not, so that there exists a morphism
$\tilde{h}:X\ra (h\mu)(K)=(h'u'\mu)(K)=(h'\mu 'u)(K)$ which
coincides with $h$ on $K$. Since $u(K)\subseteq T$ we get $(h'\mu
'u)(K)\subseteq (h'\mu ')(T)$ and the composition
$X\arr{\tilde{h}}(h'\mu
'u)(K)\hookrightarrow (h'\mu ')(T)$ is denoted by $\alpha$. If
$\hat{h}':T\ra (h'\mu ')(T)$ denotes the unique
morphism such that the composition
$T\arr{\hat{h}'}(h'u')(T)\hookrightarrow T'$
coincides with $h'\mu '$, then we have the following commutative
diagram, which can be completed with $\tilde{h'}$ due to the
universal property of pushouts:
\[\xymatrix{K\ar[r]^{\mu}\ar[d]^{u} & X\ar[d]^{u'}\ar@/^/[ddr]^{\alpha} & \\
T\ar[r]^{\mu'}\ar@/_/[drr]_{\hat{h}'} & X'\ar@{.>}[dr]^{\tilde{h'}} & \\
&& (h'\mu')(T) }
\]
But then $\tilde{h}'$ coincides with $h'$ on $T$ due to the
equality $\tilde{h}'\mu '=\hat{h}'$. That contradicts the
choice of $h'$.

Suppose next that $\mu'$ splits. Then we rewrite the pushout diagram $(*)$ as:
\[\xymatrix{0\ar[r] & K\ar[r]^{\mu}\ar[d]^{u} & X\ar[rr]^{f}\ar[d]^{\scriptsize{\left[\begin{array}{c}h\\f\end{array}\right]}} && Y\ar@{=}[d]\ar[r] & 0 \\
0\ar[r] & T\ar[r]_{\scriptsize{\left[\begin{array}{c}\id\\0\end{array}\right]}\ \ \ \ \ } & T\oplus
Y\ar[rr]_{\scriptsize{\left[\begin{array}{cc}0&\id\end{array}\right]}}
&& Y\ar[r] & 0 }
\]
We claim that $h$ satisfies the property required by assertion $2)$. Indeed, suppose not, so that  there exists $\tilde{h}:X\ra
(h\mu)(K)=u(K)$ that coincides with $h$ on $K$. Since $u$ is a
monomorphism the induced morphism
$\hat{u}:K\arr{\cong}u(K)$ is an
isomorphism, which, in addition, has the property that
$u\hat{u}^{-1}$ is the canonical inclusion $j:u(K)=(h\mu
)(K)\hookrightarrow T$. Consider the composition
\[\alpha:X\arr{\tilde{h}}u(K)\arr{\hat{u}^{-1}}K.
\]
Then $u\alpha\mu=u\hat{u}^{-1}\tilde{h}\mu=j\tilde{h}\mu=h\mu=u$.
Since $u$ is a monomorphism we get $\alpha\mu=\id_{K}$, and so
$\mu$ split, which is a contradiction.

$3)\Rightarrow 5)$ Let $f:X\ra Y$ be an epimorphism as indicated
and let $\mu:T:=\ker(f)\ra X$ be its kernel. Choose a morphism
$h:X\ra T'$ as given by assertion $3)$. Suppose that there is a
morphism $g:Y\ra T'$ such that $T+\ker(h-gf)=X$. Then we have:
$\im(h-gf)=(h-gf)(X)=(h-gf)(T)\subseteq h(T)+(gf)(T)=h(T)$.
That gives a morphism $\tilde{h}:X\ra h(T)$ such that
$i\tilde{h}=h-gf$, where $i:h(T)\ra T'$ is the canonical
inclusion. But then $i\tilde{h}\mu=(h-gf)\mu=h\mu$, which
contradicts the choice of $h$.

$5)\Rightarrow 3)$ Let $\mu:T\ra X$ be a non-split monomorphism
with  $T\in\ct$ and let $f:X\ra Y$ be its cokernel. We pick up a
morphism $h:X\ra T'\in\ct$ satisfying the hypothesis of condition $5)$. We will prove that it also satisfies the requirements of
condition $3)$. Suppose that there is a morphism $\tilde{h}:X\ra
\im(h\mu)$ such that in the diagram
\[\xymatrix{T\ar[r]^{\mu}\ar[d]_{p} & X\ar[d]^{h}\ar[dl]_{\tilde{h}}\\
\im(h\mu)\ar[r]_{i} & T' }
\]
we have $h\mu=i\tilde{h}\mu$, where $ip=h\mu$ is the
canonical epi-mono factorization of $h\mu$. Then
$(h-i\tilde{h})\mu=0$ and so there exists a morphism $g:Y\ra T'$
such that $i\tilde{h}=h-gf$. Thus
$\ker(h-gf)=\ker(i\tilde{h})\cong\ker(\tilde{h})$. We claim that
the morphism
\[\scriptsize{\left[\begin{array}{cc}\tilde{h}^k&\mu\end{array}\right]}: \ker(\tilde{h})\oplus T\ra X
\]
is an epimorphism and that will imply that $\ker(h-gf)+\ker(f)=X$,
thus yielding a contradiction. Indeed, from $\tilde{h}\mu=p$ we
get that $\tilde{h}$ is an epimorphism, and so it is the cokernel
of its kernel. Now, if $a$ is a morphism such that
$a\tilde{h}^k=0$ and $a\mu=0$, then there exists a morphism $b$
such that $b\tilde{h}=a$. Hence $0=b\tilde{h}\mu=bp$, which
implies $b=0$, and therefore $a=0$.
\end{proof}

\begin{corollary} \label{balanced wrt projective-injective}
Let $\ca$ be an abelian category and let $\ct\subseteq\ca$ be a full subcategory consisting of projective-injective objects, closed under finite direct sums and direct summands. Then $\ul{\ca}$ is balanced.
\end{corollary}
\begin{proof}
Condition $3)$ of Theorem ~\ref{characterization balanced} trivially holds.
\end{proof}
\bigskip

\section{Weak balance when $\ct$ consists of projective objects}\label{Weak balance when T consists of projective objects}
\bigskip

In this section the hypotheses on $\ct$ are the same as in section 3, but, in order to have suitable one-sided (pre)triangulated
structures on $\ul{\ca}$, further assumptions will be made in each particular result. We already know that
strong monomorphisms in $\ul{\ca}$, whenever they make sense, are
represented by epimorphisms. In order to characterize the weak
balance of $\ul{\ca}$ in case $\ct$ is functorially finite, we
will identify first those epimorphisms in $\ca$ which become
strong epis and strong monos in $\ul{\ca}$.

\begin{proposition}\label{objetos cero en right triangles2}
Let $f\in\ca(X,Y)$ be an epimorphism and suppose that $\ct$ is
covariantly finite in $\ca$. The following assertions are
equivalent:
\begin{enumerate}[1)]
\item $\ol{f}$ is a strong epimorphism in the right triangulated
category $\ul{\ca}$ \item For every (resp. some) $\ct$-preenvelope
$\mu^X:X\ra T^X$, one has that the canonical monomorphism
$\mu^X(\ker(f))\ra T^X$ splits. \item There is a decomposition
$X=T'\oplus X'$, with $T'\in\ct$, and
\[f=\scriptsize{\left[\begin{array}{cc}0 & f'\end{array}\right]}: T'\oplus X'\ra Y
\]
such that $\ca(f',T):\ca(Y,T)\arr{\sim}\ca(X',T)$ is an isomorphism for every  $T\in\ct$.
\item There is a decomposition $X=T'\oplus X'$, with $T'\in\ct$, and
\[f=\scriptsize{\left[\begin{array}{cc}0 & f'\end{array}\right]}: T'\oplus X'\ra Y
\]
such that every morphism $h:X'\ra T$ to an object of $\ct$
vanishes on $\ker(f')$.
\end{enumerate}
\end{proposition}
\begin{proof}
$1)\Leftrightarrow 2)$ Let $\mu^X:X\ra T^X$ be any
$\ct$-preenvelope. Since $f$ is an epimorphism a canonical
construction gives a cocartesian square of the form:
\[\xymatrix{X\ar[r]^{\mu^X}\ar[d]_{f} & T^{X}\ar[d]^{\pi} \\
Y\ar[r] & T^{X}/\mu^X(\ker(f))
}
\]
where $\pi$ is the canonical projection. Now, from the dual of
Proposition ~\ref{left zero}, we get that condition $1)$ holds
if and only if $T^X/\mu^X(\ker(f))\in\ct$. Since $\ct$ consists
of projective objects, the latter is equivalent to $2)$.

$3)\Leftrightarrow 4)$ Let $f':X'\ra Y$ be an epimorphism in
$\ca$. Clearly, the map $\ca (f',T):\ca (Y,T)\ra\ca (X',T)$ is an
isomorphism if and only if every homomorphism $h:X'\ra T$ vanishes
on $\ker(f')$. From that the equivalence of assertions $3)$ and
$4)$ is obvious.

$4)\Rightarrow 2)$ If we have a decomposition
$f=\scriptsize{\left[\begin{array}{cc}0&f'\end{array}\right]}:X=T'\oplus X'\ra Y$ as
given by condition $4)$, then $f'$ is necessarily an epimorphism in $\ca$. Let now
$\mu^X=\begin{pmatrix}u & v
\end{pmatrix}:T'\oplus X'\ra T^X$ be any $\ct$-preenvelope. One readily sees
that $u$ is a section, so that $\mu^X$ can be written as a matrix
$\mu^X=\scriptsize{\left[\begin{array}{cc}\id&\varphi\\0&\psi\end{array}\right]}:T'\oplus
X'\ra T'\oplus T''=T^X$. One then checks that
$\psi:X'\ra T''$ is a $\ct$-preenvelope and, in
particular, we have $\varphi =\rho\psi$, for some morphism $\rho
:T''\ra T'$. Now we have
$\mu^X(\ker(f))=\scriptsize{\left[\begin{array}{cc}\id&\rho\psi\\0&\psi\end{array}\right]}(T'\oplus
\ker(f'))$. By hypothesis, $\psi (\ker(f'))=0$, from which it
follows that $\mu^X(\ker(f))=T'\oplus 0$ is a direct summand of
$T^X$.

$2)\Rightarrow 4)$ Let us fix a $\ct$-preenvelope $\mu^X:X\ra T^X$ satisfying condition $2)$. Then the composition of inclusions
$\mu^X(\ker(f))\hookrightarrow \im(\mu^X)\hookrightarrow T^X$ is a section in $\ct$, which implies that $T'=:\mu^X(\ker(f))$ is also a
direct summand of $\im(\mu^X)$, so that we have a decomposition $\im(\mu^X)=T'\oplus V$. Looking now at the induced epimorphisms
$\tilde{\mu}:X\twoheadrightarrow \im(\mu^X)=T'\oplus V$ and $\ker(f)\twoheadrightarrow \mu^X(\ker(f))=T'$,  and bearing in mind
that $T'$ is a projective object, we can assume without loss of generality that there are decompositions $X=T'\oplus X'$ and
$f=\scriptsize{\left[\begin{array}{cc}0 & f' \end{array}\right]}:T'\oplus X'\ra Y$ (\ie $T'\subseteq \ker(f)$). Then $\mu^X$ gets identified with
$\scriptsize{\left[\begin{array}{cc}1 & 0\\ 0 & \mu^{X'}\end{array}\right]}:T'\oplus X'\ra T'\oplus T^{X'}=T^{X}$, where $\mu^{X'}:X'\ra T^{X'}$ is a $\ct$-preenvelope such that $\mu^{X'}(\ker(f'))=0$. Since every morphism $h :X'\ra T$, with $T\in\ct$, factors through $\mu^{X'}$
assertion $4)$ follows.
\end{proof}

\begin{proposition} \label{epimorphisms which become strong monos}
Let $f\in\ca (X,Y)$ be an epimorphism and suppose that $\ct$ is contravariantly finite in
$\ca$. In the left triangulated category $\ul{\ca}$, the morphism $\ol{f}$ is a strong
monomorphism if and only if $\ker(f)\in\ct$.
\end{proposition}
\begin{proof}
If $p_Y:T_Y\ra Y$ is any $\ct$-precover, then the pullback of $f$ and $p_{Y}$
\[\xymatrix{Z\ar[r]^{g}\ar[d]_{r} & T_Y\ar[d]^{p_Y} \\
 X\ar[r]_{f} & Y
}
\]
gives rise to a  short exact sequence $0\ra \ker(f)\ra Z\arr{g}
T_Y\rightarrow 0$ in $\ca$, which splits due to the fact that
$T_Y$ is projective. Then $\ker(f)$ belongs to $\ct$ if and only if $Z$ belongs to $\ct$. Then the result follows from Proposition
~\ref{left zero}.
\end{proof}

\begin{corollary} \label{epimorphisms with zero left-and-right vertex}
Suppose that $\ct$ is functorially finite in $\ca$ and let $f:X\ra Y$ be an epimorphism
in $\ca$. The following assertions are equivalent:

\begin{enumerate}[1)]
 \item $\ol{f}$ is a strong monomorphism and a strong epimorphism in the pretriangulated category $\ul{\ca}$.
\item $\ker(f)\in\ct$ and, for some $\ct$-preenvelope $\mu^X:X\ra T^X$, one has that $\mu^X(\ker(f))$ is a direct summand of $T^X$.
\item $\ker(f)\in\ct$ and there is a decomposition
$f=\scriptsize{\left[\begin{array}{cc}0&f'\end{array}\right]}:X=T'\oplus X'\ra Y$ such
that every morphism $h :X'\ra T$ to an object of $\ct$ vanishes on $\ker(f')$
\end{enumerate}

\end{corollary}
\begin{proof}
Direct consequence of Propositions \ref{objetos cero en right triangles2} and
\ref{epimorphisms which become strong monos}.
\end{proof}

All these preliminaries lead to the following main result of this
section.

\begin{theorem} \label{characterization weakly-balanced}
Let $\ca$ be an abelian category and $\ct\subseteq\ca$ be a
functorially finite full subcategory consisting of projective
objects and which is closed under taking direct summands. Consider
the following assertions:
\begin{enumerate}[1)]
\item For every object $T\in\ct$, there is a monomorphism
$T\rightarrowtail E$, where $E$ is an injective(-projective)
object of $\ca$ which belongs to $\ct$. \item For every
$T\in\ct\setminus\{0\}$, there is a nonzero morphism $\varphi:T\ra
T'$, with $T'\in\ct$, which factors through an injective object of
$\ca$. \item The pretriangulated category $\ul{\ca}=\ca/<\ct>$ is
weakly balanced. \item If $f:X\ra Y$ is an epimorphism in $\ca$
such that $\ker(f)\in\ct$ and every morphism $h:X\ra T\in\ct$
vanishes on $\ker(f)$, then $f$ is an isomorphism in $\ca$. \item
If $0\neq T\arr{j}X$ is a monomorphism in $\ca$, with $T\in\ct$,
then there is a morphism $h :X\ra T'$ such that $h j\neq 0$,  for
some $T'\in\ct$. \item If $Y$ is an object of $\ca$ such that the
covariant functor
\[\Ext_{\ca}^1(Y,?)_{|\ct}:\ct\ra \Mod\Z
\]
is nonzero, then $\Ext_{\ca}^1(Y,?)_{|\ct}$ does not contain a nonzero representable subfunctor.
\end{enumerate}

Then $1)\Rightarrow 2)\Rightarrow 3)\Leftrightarrow 4)\Leftrightarrow 5)\Leftrightarrow (6)$ and, in case $\ca$ has enough injecives, all assertions $2)-6)$ are equivalent.
\end{theorem}
\begin{proof}
$1)\Rightarrow 2)$ is clear.

$2)\Rightarrow 4)$ Let $f:X\ra Y$ be an epimorphism in $\ca$
satisfying the hypothesis of $(4)$. Suppose that $f$ is not an
isomorphism. Then $\ker(f)\in\ct\setminus\{0\}$, and by $(2)$
there exists a non-zero morphism $\varphi:\ker(f)\ra T'$, for some
$T'\in\ct$, which factors through an injective object of $\ca$.
Then $\varphi$ extends to $X$ and this contradicts the hypothesis
of $(4)$.

$3)\Leftrightarrow 4)$ is a direct consequence of Corollary
\ref{epimorphisms with zero left-and-right vertex} and Proposition
\ref{epimorphisms which become stable-isos}, bearing in mind that
strong epimorphisms of $\ul{\ca}$ can be represented by
epimorphisms of $\ca$.

$4)\Leftrightarrow 5)$ is clear.

$4)\Rightarrow 6)$ Let $Y\in\ca$ be an object such that
$\Ext_{\ca}^1(Y,?)_{|\ct}:\ct\ra\Mod\Z$ is a nonzero functor.
Suppose that we have a monomorphism of functors $\mu :\ct (T,?)\ra
\Ext_{\ca}^1(Y,?)_{|\ct}$.  We choose the element
$\tilde{\mu}\in \Ext_{\ca}^1(Y,T)$ corresponding to $\mu$ by
Yoneda's lemma. Then $\tilde{\mu}$ represents a short exact
sequence $0\ra T\arr{j}X\arr{f}Y\ra 0$ in $\ca$. Then $f$ is an
epimorphism in $\ca$ such that $\ker(f)\in\ct$. On the other hand,
if $h :X\ra T'$ is a morphism to an object of $\ct$, then
$\mu_{T'}:\ct (T,T')\ra \Ext_{\ca}^1(Y,T')$ maps $hj$
onto the (exact) lower row of the following pushout diagram:
\[\xymatrix{0\ar[r] & T\ar[r]^{j}\ar[d]_{hj} & X\ar[r]\ar[d] & Y\ar[r]\ar@{=}[d] & 0 \\
0\ar[r] & T'\ar[r] & X'\ar[r] & Y\ar[r] & 0
}
\]
That lower row splits because $hj$ factors through $j$.
Therefore $\mu_{T'}(hj)=0$. The monomorphic condition of
$\mu$ implies that $hj=0$. Thus, every morphism
$h:X\ra T'$, with $T'\in\ct$, vanishes on $T=\ker(f)$.
Now assetion $(4)$ says that $f$ is an isomorphism in $\ca$ and,
hence, that $T=0$.

$6)\Rightarrow 4)$ Let $f:X\ra Y$ be an epimorphism as indicated
in assertion $4)$ and consider the associated short exact
sequence $0\ra T\arr{j} X\arr{f}Y\ra 0$ in $\ca$. By Yoneda's
lemma, we get an induced morphism of functors $\mu :\ct (T,?)\ra
\Ext_{\ca}^1(Y,?)_{|\ct}$ (functors $\ct\ra\Mod\Z$). We claim
that $\mu$ is a monomorphism and then the hypothesis will give
$\ct (T,?)=0$, whence $T=0$, which will end the proof of this
implication. Let us prove our claim. If $T'\in\ct$ is any object
then $\mu_{T'}:\ct (T,T')\ra \Ext_{\ca}^1(Y,T')$ maps a
morphism $\varphi:T\ra T'$ onto the (exact) lower row of the
following pushout diagram:
\[\xymatrix{0\ar[r] & T\ar[r]^{j}\ar[d]_{\varphi} & X\ar[r]\ar[d] & Y\ar[r]\ar@{=}[d] & 0 \\
0\ar[r] & T'\ar[r] & X'\ar[r] & Y\ar[r] & 0
}
\]
One has that $\mu_{T'}(\varphi )=0$ if and only if  that lower
row splits. But this is equivalent to say that $\varphi$ is the
restriction to $T=\ker(f)$ of a morphism $h :X\ra T'$. In that
case, the hypothesis on $f$ says that $\varphi =0$. Therefore
$\mu$ is a monomorphism of functors as desired.

We finally prove $5)\Rightarrow 2)$ in case $\ca$ has enough
injectives. Let us take $0\neq T\in\ct$.  There is a monomorphism
$j:T\ra E$, with $E$ injective. Now assertion $5)$ says that
there exists a morphism $h:E\ra T'$, with $T'\in\ct$, such that
$\varphi =hj\neq 0$. That ends the proof.
\end{proof}

%%\begin{remark}
%%The reader will have notice a certain lack of rigour in assertion
%7 above. Namely, it may happen that the extensions of two objects
%do not form a set. The usual 'increasing of universe', passing
%from the category $Ab$ of (small) abelian groups to the category
%of 'large abelian groups', solves the problem. In order to avoid
%these set-theoretical technicalities we have chosen to state the
%theorem as though $Ext_\ca(-,-)$ takes values in $Ab$, which is
%the most usual case in practice.
%\end{remark}

We shall end the section by characterizing (weak) balance in
abelian categories on which subobjects of projective objects are
projective. We denote by $Sub(\ct )$ the full subcategory having
as objects those which are subobjects of objects in $\ct$. The
following is an auxiliary lemma.

\begin{lemma} \label{covariantly-finite in hereditary}
Let $\ch$ be an abelian category in which subobjects of projective
objects are projective and let $\ct$ be a full subcategory of
$\ch$ consisting of projective objects and closed under direct
summands. The following assertions are equivalent:

\begin{enumerate}
\item $\ct$ is covariantly finite in $\ch$ \item $\op{Sub}(\ct )$ is
covariantly finite in $\ch$ and every $P\in Sub(\ct )$ has a
$\ct$-preenvelope  \item $(^\perp\ct ,\op{Sub}(\ct ))$ is a (split)
torsion theory in $\ch$ and every $P\in \op{Sub}(\ct )$ has a
$\ct$-preenvelope
\end{enumerate}
\end{lemma}
\begin{proof}
$1)\Longrightarrow 2)$ If $\mu^X:X\ra T^X$ is a
$\ct$-preenvelope and we take its epi-mono factorization
$X\arr{p}P\stackrel{j}{\hookrightarrow}X$,
then $p$ is a $\op{Sub}(\ct )$-preenvelope.

$2)\Rightarrow 1)$ If $\lambda :X\ra P$ and $\mu
:P\ra T^P$ are a $\op{Sub}(\ct )$-preenvelope and a
$\ct$-preenvelope, respectively, then $\mu\lambda$ is a
$\ct$-preenvelope.

$2)\Leftrightarrow 3)$ is a direct consequence of the dual
of Lemma \ref{(co- contra-)variantly finite Serre classes}
\end{proof}

Recall that a torsion theory $(\cu ,\cv )$ in $\ca$ is called \emph{hereditary}  when $\cu$  is closed under taking subobjects.

\begin{proposition} \label{hereditary weakly balanced}
Let $\ch$ be an abelian category in which subobjects of projective
objects are projective, and let $\ct\subseteq\ch$ be a covariantly
finite full subcategory of $\ch$ consisting of projective objects
and closed under taking direct summands. The following assertions
are equivalent:

\begin{enumerate}[1)]
\item $\ul{\ch}=\ch/<\ct>$ is  balanced. \item  $^\perp\ct$ is
closed under taking subobjects. \item The pair
$(^\perp\ct,\op{Sub}(\ct))$ is a hereditary (split) torsion theory
in $\ch$.
\end{enumerate}

When $\ct$ is contravariantly finite in $\ch$, the above
assertions are equivalent to:

\begin{enumerate}[4)]
\item  $\ul{\ch}$ is weakly balanced
\end{enumerate}

When $\ch$ has enough injectives,  assertions $1), 2)$ and $3)$ are
also  equivalent to:

\begin{enumerate}[5)]
\item For every object $T\in\ct$, there is a monomorphism
$T\rightarrowtail E$, where $E$ is an injective(-projective)
object of $\ch$ which belongs to $\ct$.
\end{enumerate}
\end{proposition}

\begin{proof}
We start by proving that, under anyone  of conditions $1)$, $3)$ or $4)$,  if
\[\lambda:=\scriptsize{\left[\begin{array}{c}u\\ v\end{array}\right]}:T\ra Y\oplus P
\]
is a monomorphism, where $T\in\ct$, $Y\in ^\perp\ct$ and $P$ is
projective, then $v$ is also a monomorphism. Indeed, since
subobjects of projective objects are projective we get that
$T_2:=\im(v)$ is projective and the induced epimorphism
$\tilde{v}:T\ra T_2$ splits. Then we can decompose $T\cong
T_1\oplus T_2$, where $T_1:=\ker(v)$, and then rewrite $\lambda$
as $\scriptsize{\left[\begin{array}{cc}u_1 & u_2
\\ 0 & v_2
\end{array}\right]}$. Since $\lambda$ is a monomorphism it follows
easily that $u_1$ is a monomorphism. Under conditions $1)$ or
$3)$, using condition $3)$ of Theorem ~\ref{characterization
balanced} or condition $5)$ of Theorem \ref{characterization
weakly-balanced},  one readily sees that necessarily $T_1=0$.
Under condition $2)$, one has that $T_1\in \ct\cap ^\perp\ct$
since the torsion theory $(^\perp\ct,\op{Sub}(\ct))$ is
hereditary. But then $T_1=0$ as well.

$1)\Rightarrow 2)$ By Lemma \ref{covariantly-finite in
hereditary}, we know that $(^\perp\ct, \op{Sub}(\ct ))$ is a split
torsion theory in $\ch$. Let $i:X\ra Y$ be a monomorphism, where
$Y\in ^\perp\ct$, and take any morphism $h:X\ra T$, with
$T\in\ct$. We then take the bicartesian square:
\[\xymatrix{X\ar[r]^{i}\ar[d]_{h} & Y\ar[d]^{\eta} \\
  T\ar[r]_{\lambda} & V
}
\]
Then $\lambda$ is a monomorphism. We decompose $V=V'\oplus P$,
where $V'\in ^\perp\ct$ and $P\in \op{Sub}(\ct )$. According to the
first paragraph of this proof, we can rewrite
$\lambda=\scriptsize{\left[\begin{array}{c}\gamma \\ j\end{array}\right]}:T\ra V'\oplus P$,
 where
$j:T\ra P$ is a monomorphism. On the other hand, since $(^\perp\ct, \op{Sub}(\ct ))$
is a torsion theory, we have $\ca (Y,P)=0$. Then
$\eta=\scriptsize{\left[\begin{array}{c}u\\ 0 \end{array}\right]}$.
The commutativity of the square implies that $jh=0$ and, since $j$ is a
monomorphism, we conclude that $h=0$. Therefore $X\in ^\perp\ct$.

$2)\Rightarrow 3)$ is a direct consequence of Lemma \ref{covariantly-finite in hereditary}.

$3)\Rightarrow 1)$ We shall check assertion $3)$ of Theorem
~\ref{characterization balanced}. Let $\mu :T\rightarrowtail X$ be
a non-split monomorphism in $\ch$ with $0\neq T\in\ct$.  By
hypothesis, we can decompose $X=X'\oplus P$, where $X'\in
^\perp\ct$ and $P\in \op{Sub}(\ct )$. Then the first paragraph of
this proof says that
we have $\mu=\scriptsize{\left[\begin{array}{c}u \\
v\end{array}\right]}$, where $v:T\ra P$ is a monomorphism. But,
since $P\in \op{Sub}(\ct )$, we have a monomorphism $w:P\ra T'$,
with $T'\in\ct$. Clearly $wv\neq 0$ (it is actually a
monomorphism). We then claim that
$h:=\scriptsize{\left[\begin{array}{cc}0&w\end{array}\right]}:X'\oplus
P=X\ra T'$ is a morphism such that no morphism
$\tilde{h}:X=X'\oplus P\ra (h\mu )(T)=wv(T)$ coincides
with $h$ on $T$. To see that notice that, if it exists, such a
$\tilde{h}$ can be written as $\tilde{h}=\begin{pmatrix}0 & \beta
\end{pmatrix}:X'\oplus P\ra wv(T)$, because $X'\in
^\perp\ct$ and $wv(T)\cong T\in\ct$. If we denote by $i$ the
canonical inclusion $wv(T)\hookrightarrow T'$ and by
$\xi:T\arr{\cong}wv(T)=h\mu (T)$ the
isomorphism induced by $h\mu$ (or $wv$), then we have $i\xi =h\mu=i\tilde{h}\mu =i\scriptsize{\left[\begin{array}{cc}0 & \beta\end{array}\right]}
\scriptsize{\left[\begin{array}{c}u\\v \end{array}\right]}=i\beta v$. Since $i$ is a
monomorphism we conclude that $\xi =\beta v$, so that $v$ is a
section. But then $\mu=\scriptsize{\left[\begin{array}{c}u\\v\end{array}\right]}$ is also
a section, which contradicts the hypothesis. Therefore $\tilde{h}$
cannot exists and the proof of this implication is finished.

Suppose now that $\ct$ is contravariantly finite in $\ch$. Then
$1)\Rightarrow 4)$ clearly holds. The proof of $4)\Rightarrow 2)$ is identical to that of $1)\Rightarrow 2)$,  done above.

Finally, suppose  that $\ca$ has enough injectives.

$5)\Rightarrow 1)$ We check condition $3)$ of Theorem
\ref{characterization balanced}. Let $\mu :T\rightarrowtail X$ be
a non-split monomorphism. By hypothesis, there is another
monomorphism $j:T\rightarrowtail E$, where $E$ is an injective
object of $\ch$ belonging to $\ct$. For simplicity, we view $j$ as
an inclusion. Now $j$ factors through $\mu$, so that we get a
morphism $h:X\ra E$ such that $h\mu =j$. If there were
a morphism $\tilde{h}:X\ra h\mu (T)=j(T)=T$ agreing
with $h$ on $T$, then $\tilde{h}$ would be a retraction for $\mu$,
which is absurd.

$3)\Rightarrow 5)$ Let $0\neq T\in\ct$ be any nonzero object. We
then have a monomorphism $\lambda :T\ra E$, where $E$ is an
injective object of $\ca$. We decompose $E=E'\oplus P$, with
$E'\in ^\perp\ct$ and $P\in \op{Sub}(\ct )$. Notice that then $P$
is injective(-projective) and necessarily belongs to $\ct$. If we
write now $\lambda =\scriptsize{\left[\begin{array}{c}u \\
v\end{array}\right]}:T\ra E'\oplus P$, then the first paragraph of
this proof says that $v$ is a monomorphism. That ends the proof.
\end{proof}

\begin{remark}
On an arbitrary abelian category $\ca$,  Yoneda defined the big
abelian group  $\Ext_{\ca}^n(X,Y)$ of  $n$-extensions (see
\cite[III.5]{MacLane1975}), for all objects $X,Y\in\ca$, which was
functorial on both variables.  Due to Verdier \cite[Ch. III, 3.2.12]{Verdier1996},
$\Ext_\ca^n(-,-)=\Hom_{\cd\ca}(X,Y[n])$, where $\cd\ca$ is the (possibly
large) derived category of $\ca$. As a consequence, the long exact
sequences of $\Ext$ hold in every abelian category. In particular,
if $\ca =\ch$ is hereditary (\ie $\Ext_{\ch}^2(-,-)=0$), then
every subobject of a projective object in $\ch$ is projective.
Thus our Proposition \ref{hereditary weakly balanced} applies to
hereditary abelian categories.
\end{remark}
\bigskip

\section{Some examples}\label{Some examples}
\bigskip

In this section we illustrate the previous results with some
examples in module categories. Unless said otherwise, all rings
are associative with unit and modules are left modules. For a
given ring $R$, we denote by $R\Mod$ (resp. $R\mod$) the category
of all (resp. finitely presented) left $R$-modules.

\begin{example}
Let $R$ be a left semihereditary right coherent ring. Let $\ct
=R\op{proj}$ the full subcategory of the abelian category
$\ch=R\mod$ consisting of finitely generated projective modules.
By \cite[Corollary 3.11]{DingChen1993}, we know that $\ct$ is
covariantly finite in $\ch$. On the other hand, the left
semihereditary condition of $R$ gives that in $\ch$ subobjects of
projective objects are projective, so that we are in the situation
of Proposition \ref{hereditary weakly balanced}. Then
$R\ul{\mod}={R\mod}/{R\op{proj}}$ is (weakly) balanced if and only
if the class $^\perp\ct =\{M\in R\mod\mid M^*=\Hom_R(M,R)=0\}$ is
closed under (finitely generated) submodules.

Two particular cases of this situation are those in which $R$ is a commutative principal ideal domain or $R=k\vec{A_{n}}$ is the path algebra
of the Dynkin quiver $\vec{A_{n}}:$ $1\ra 2\ra ...\ra n-1\ra n$ over an algebraically closed field $k$.
\end{example}

It is well-known that the above path algebra $k\vec{A_{n}}$ is isomorphic to $T_n(k)$, where, for a ring $R$ and a natural number $n$, $T_n(R)$
denotes the ring of (upper) triangular $n\times n$ matrices with coefficients in $R$.  When one moves from finitely presented to
arbitrary modules, only this example survives.

\begin{example}
Let $H$ be a left hereditary ring. The following assertions are equivalent:
\begin{enumerate}
\item $H\op{Proj}$ is covariantly finite in $H\Mod$ and
$H\ul{\Mod}=H\Mod/H\op{Proj}$ is a (weakly) balanced category.
\item $H$ is Morita equivalent to a finite direct product
$T_{n_1}(D_1)\times ...\times T_{n_r}(D_r)$, for some division
rings $D_i$ and natural numbers $n_i$
\end{enumerate}
\end{example}
\begin{proof}
$1)\Rightarrow 2)$ Since $H\op{Proj}$ is covariantly finite in
$H\Mod$ it is closed under taking products, so that $H$ is
semiprimary (cf. \cite{Chase1960} and \cite{Small1967}). On the
other hand,  since assertion $4)$ of Proposition \ref{hereditary
weakly balanced} holds, we get that $E(_HH)$ is projective, so
that $H$ is left QF-3 in the sense of Ringel and Tachikawa (see
\cite[Proposition 4.1]{Tachikawa1973} ). Then, by
\cite[Theorem and Remark 1]{ColbyRutter1968}, one concludes that
$R$ is Morita equivalent to the indicated direct product of matrix
rings.

$2)\Rightarrow 1)$ There is no loss of generality in
assuming  that $H=T_n(D)$, for some natural number $n$ and some
division ring $D$. Since $H$ is  (two-sided) Artinian, the class
$H\op{Proj}$  is covariantly finite in $H\Mod$ (see
\cite[Corollary 3.5]{RadaSaorin1998}. On the other hand,
 $H$ is a serial ring and the column  module $E=\scriptsize{\left[\begin{array}{c} D\\ D\\ \cdot \\ \cdot \\D\end{array}\right]}=D^n$ is the unique injective-projective
indecomposable left $H$-module, which cogenerates all other
projective left $H$-modules. Then assertion $5)$ of Proposition
\ref{hereditary weakly balanced} holds, so that $H\ul{\Mod}$ is
balanced.
\end{proof}

In the situation of Proposition \ref{hereditary weakly balanced}, whenever $\ct$ is closed for subobjects, one has an equivalence of categories $\ul{\ch}=\frac{\ch}{<\ct >}\cong ^\perp\ct$, so that $\ul{\ch}$ is an abelian category. That is no longer true if $\ct$ is not closed for subobjects, as the following example shows.

\begin{example} \label{hereditary non-abelian non-triangulated}
Let $H=k\vec{A_{3}}$,  where $\vec{A_{3}}$ is the Dynkin quiver $1\ra 2\ra 3$ and $k$ is an algebraically closed field. Let $e_i$ be the primitive idempotent of $H$ given by the vertex $i$ ($i=1,2,3$) and take $\ch =H\mod$ and $\ct =\op{add}(He_3)$. Then $\ul{\ch}={H\mod}/{<\ct >}$ is a balanced category which is neither abelian nor triangulated (for its usual pretriangulated structure).
\begin{proof}
Due to Corollary \ref{balanced wrt projective-injective},
$\ul{\ch}$ is balanced. If $f:He_2\twoheadrightarrow
S_2=\frac{He_2}{J(H)e_2}$ is the canonical projection, then, by
Proposition \ref{epimorphisms which become stable-epis}, we have
that $\ol{f}$ is not an epimorphism in $\ul{\ch}$. If this latter
category were abelian, we would have a morphism $g:S_2\ra X$ in
$\ch$ such that $0\neq\ol{g}=\ol{f}^c$. But  every such morphism
$g$ necessarily factors through the (unique up to multiplication
by scalar) irreducible morphism starting at $S_2$, which is the
injective envelope map $j:S_2\hookrightarrow
E(S_2)=\frac{He_3}{Soc(He_3)}$. Since $jf$ factors through $He_3$,
we have $\ol{j}\ol{f}=0$ and from that one derives that, in case
of existence of $\ol{f}^c$, one necessarily has $\ol{f}^c=\ol{j}$.
But this is absurd for $\ol{j}$ is not an epimorphism in
$\ul{\ch}$: if
$p:E(S_2)=\frac{He_3}{Soc(He_3)}\twoheadrightarrow\frac{He_3}{J(H)e_3}=S_3$
is the canonical projection, then $pj=0$ but $\ol{p}\neq 0$.

Notice also that $\ct^\perp\neq 0$ (e.g. $He_1\in\ct^\perp$), so that the loop functor $\Omega$, given by the pretriangulated
structure of $\ul{\ch}$, vanishes on some nonzero objects. Therefore $\ul{\ch}$ is not a triangulated category either.
\end{proof}

\end{example}

Corollary \ref{balanced in Serre case} and Proposition \ref{hereditary weakly balanced} might induce the reader to believe that `weak balance' (when it makes sense) and `balance' are synonimous concepts for stable categories. Our final example shows that it is not the case, even when $\ct$ consists of projective objects.

\begin{example} \label{weakly balanced not balanced}
Let $k$ be an algebraically closed field and $A$ be the (finite dimensional) $k$-algebra given by the quiver
\[\xymatrix{1\ar[d]_{\delta} & 2\ar[l]_{\alpha}\ar[r]^{x} & 5\ar[r]^{y} & 6 \\
                  4\ar[r]_{\gamma} & 3\ar[u]_{\beta} &&
}
\]
with the following set of monomial relations  $\rho =\{xy,\beta
x,\gamma\beta ,\alpha\delta\gamma ,\beta\alpha\delta\}$. Then
$A\op{proj}$ is functorially finite in $A\mod$. We claim that
$A\ul{\mod}={A\mod}/{A\op{proj}}$ is weakly balanced but  not
balanced. To see that, first notice that assertion $2)$ of Theorem
\ref{characterization weakly-balanced} is equivalent in this case
to the property that, for every indecomposable projective left
$A$-module $P$, the canonical restriction  map $\Hom_A(E(P),A)\ra
\Hom_A(P,A)$ has a nonzero image. This property is trivially
satisfied whenever $E(P)$ is projective. In our situation, that is
the case for $P=Ae_i$, with $i\in\{1,2,3,6\}$. On the other hand,
we have $E(Ae_4)\cong E(Ae_5)\cong \frac{Ae_4\oplus Ae_5}{N}=:E$,
where $N$ is the cyclic submodule of $Ae_4\oplus Ae_5$ generated
by $(-\alpha\delta ,x)$. One readily sees that $\Hom_A(E,A)$ is a
$2$-dimensional vector space generated by the morphisms
$f:\ol{(a,b)}\mapsto a\gamma$ and $g:\ol{(a,b)}\mapsto by$. If now
$i:Ae_4\hookrightarrow E$ ($a\mapsto\ol{(a,0)}$) and
$j:Ae_5\hookrightarrow E$ ($b\mapsto\ol{(0,b)}$) are the canonical
inclusions, then clearly $fi\neq 0$ and $gj\neq 0$. That proves
that $A\ul{\mod}$ is weakly balanced.

We shall see now that condition $3)$ of Theorem \ref{characterization balanced} is not satisfied by the canonical inclusion $\mu =i:Ae_4\hookrightarrow E$. Indeed, by the above paragraph, any morphism $h =E\ra _AA$ is of the form $h =k_1g+k_2h$, for some $k_1,k_2\in k$. Then it maps
$\ol{(a,b)}\mapsto k_1a\gamma+k_2by$ and, hence, $(h\mu)(a)=k_1a\gamma$. In case $k_1=0$, clearly $h$ does not satisfy
condition $3)$ of Theorem \ref{characterization balanced}. But if
$k_1\neq 0$ then $\im(h\mu)=\im(f)$ and, choosing
$\tilde{h}=k_1f:E\ra \im(f)$,
$\ol{(a,b)}\mapsto k_1a\gamma$,  one clearly has that
$\tilde{h}$ coincides with $h$ on $Ae_4$. That proves that
$A\ul{\mod}$ is not balanced.

\end{example}

\end{document}